\documentclass[12pt,leqno]{amsart}
\usepackage{amssymb,amsmath,amsfonts,amsthm,latexsym,url,enumerate}

\begin{document}

\newtheorem{theorem}{Theorem} %
\newtheorem{lemma}[theorem]{Lemma} %
\newtheorem{claim}[theorem]{Claim} %
\newtheorem{corollary}[theorem]{Corollary} %
\newtheorem{proposition}[theorem]{Proposition}
\newtheorem{question}[theorem]{Question} %
\newtheorem{definition}[theorem]{Definition} %

\theoremstyle{remark} %
\newtheorem*{remark}{Remark} %

\def\C{{\mathbb C}}
\def\F{{\mathbb F}}
\def\N{{\mathbb N}}
\def\Z{{\mathbb Z}}
\def\Q{{\mathbb Q}}
\def\Int{{\operatorname{Int}}}

\def\({\left(}
\def\){\right)}


\title{\bf Power Moments of Kloosterman Sums}


\author{Ke Gong}
\address{Department of Mathematics, Henan University, Kaifeng 475004, P.R. China}
\email{kg@henu.edu.cn}
\thanks{}

\author{Willem Veys}
\address{KU Leuven, Dept. Wiskunde, Celestijnenlaan 200B, 3001 Leuven, Belgium}
\email{wim.veys@wis.kuleuven.be}
\thanks{}

\author{Daqing Wan}
\address{Department of Mathematics, University of California, Irvine, CA92697-3875, USA}
\email{dwan@math.uci.edu}
\thanks{}

\thanks{The first author is supported by the NSFC grant 11201117. The second author is supported  by the Research Project G.0939.13N of the Research Foundation - Flanders (FWO). The third author is partially supported by the Simons Fellowship 304279 and by the NSF grant CCF-1405564.
We would like to thank Wouter Castryck and Joshua Hill for doing computer experiments, which were quite helpful
in obtaining the correct formulas in our Igusa zeta function computation.}

\date{}
\maketitle

\begin{abstract}In this paper we  give an essential  treatment for power moments of Kloosterman sums over the residue class ring $\Z/q\Z$ of $q$ elements.
When $q$ is a large enough power of a prime,
we prove concrete formulas using computations with Igusa zeta functions.  As a consequence, we show that there are explicit formulas for power moments of Kloosterman sums whenever $q$ is powerful.
\end{abstract}


\section{Introduction}


Let $q\ge 3$ be a positive integer and $\zeta=e^{2\pi i/q}$ the $q$-th root of unity. For arbitrary integers $u$ and $v$, the classical Kloosterman sums are defined by
$$
K(u,v;q) = \sum_{\substack{0<x\le q \\ (x,\,q)=1}}\zeta^{ux+v\bar x},
$$
where $\bar x$ denotes the multiplicative inverse of $x$ in $(\Z/q\Z)^*$.

The following fundamental properties of $K(u,v;q)$ are well-known, see Sali\'e~\cite[\S 1]{Salie} and Esterman~\cite[p.\,91]{Esterman}:
\begin{enumerate}[(1)]
 \item $K(u,v;q)$ is real;
 \item $K(u,v;q)=K(v,u;q)$;
 \item $K(u,v;q)=K(1,uv;q)$ when $(u,q)=1$ and $v$ is arbitrary;
 \item given integers $v$, $q_1$, $q_2$ with $(q_1,q_2)=1$, there exist integers $v_1$ and $v_2$ such that $v\equiv v_1q_1^2+v_2q_2^2\pmod q$, and $$K(u,v;q_1q_2)=K(u,v_1;q_1)K(u,v_2;q_2).$$
\end{enumerate}

The study of Kloosterman sums $K(u,v;q)$ is reduced by (4) to the prime-power case $K(u,v;p^m)$, where $m\ge 1$ is an integer and $p$ is a prime. Indeed,  the following is known at least to Esterman~\cite[p.\,91]{Esterman}.

Let $q=p_1^{m_1}p_2^{m_2}\cdots p_s^{m_s}$ with $p_1,p_2,\ldots,p_s$ distinct primes. Given $v$, we can determine $v_1,v_2,\ldots,v_s$ such that for all $u$ we have
$$
K(u,v;q)=\prod_{i=1}^s K(u,v_i;p_i^{m_i}).
$$

For a fixed $v$, the individual sum $K(1,v; q)$ has been studied extensively in the literature.
An important problem is to understand the distribution of the Kloosterman sum $K(1, v; q)$ as $v$ varies.
This problem is equivalent to understanding the $n$-th power moment of the Kloosterman sums defined by
$$
S_n(q)=\sum_{u=1}^q K(u,1;q)^n = \sum_{v=1}^q K(1,v;q)^n,
$$
where $n$ is a positive integer.  Ideally, one would like to have an explicit formula for the $n$-th moment $S_n(q)$.
This is indeed the case for small $n$, as one checks directly that
$$S_1(q)=0, \ S_2(q) = \varphi(q)q,$$
where $\varphi(q)$ is the Euler function.
As $n$ grows, one expects that the moments $S_n(q)$ become increasingly complicated. This is indeed the case for general $q$.
A surprising consequence of our work is that there are also explicit formulas for $S_n(q)$  for higher $n$ if $q$ is powerful compared to $\log_p n$.
We state two results in this direction. For a prime $p$, we let $v_p(q)$ denote the $p$-adic valuation of $q$.

\begin{theorem} Let $n\geq 1$ be odd. If $q$ has a prime factor $p$ such that $v_p(q)  >  \log_p n +1$, then
$S_n(q)=0$.
\end{theorem}

\begin{theorem} Let $n\geq 2 $ be even. If $q$ is odd and each prime factor $p$ of $q$ satisfies
$v_p(q)  >  \log_p \frac{n}{2} +1$, then
$$
S_n(q) = \prod_{p|q} \binom {n-1} {\frac n2 -1} \frac{p-1}p p^{(\frac n2 +1) v_p(q)}.
$$
\end{theorem}
\noindent
If $q$ is even, there is a similar result that we do not state here for simplicity.

\medskip
Here we take the point of view of fixing $n$ and letting $q$ vary.  This yields uniform formulas for sufficiently powerful $q$.  Of course, one could also think about fixing $q$ and letting $n$ vary; in that case our formulas are valid for finitely many $n$.

\medskip
{\noindent} {\sl Remark.}  Motivated by applications in cryptography and coding theory \cite{Carlet}, it is a challenging problem to give a sharp
lower bound for the quantity
$$M_q = \max_{0\leq u\leq q-1}|K(u, 1; q)|.$$
To get a good clean lower bound for $M_q$, it is essential to get a clean explicit formula for the higher moments $S_n(q)$.
(Indeed, by the definition of $S_n(q)$, one gets the inequality $S_n(q) \leq q M_q^n$.)
Our result fits exactly this purpose when $q$ is powerful.  We hope to study this type of applications in a more systematic way
in a future paper.

\medskip

The study of the power moment $S_n(q)$ can be easily reduced to the case when $q$ is a prime power.
For any integer $q\ge 3$ and $v$ with $(v,q)=1$, let $q=\prod_{i=1}^s p_i^{m_i}$, and
$$
u\equiv u_i\pmod{p_i^{m_i}},\qquad 1\le i\le s.
$$
It is clear that when $u_i$ runs through a complete residue system modulo $p_i^{m_i}$, $u$ runs through a complete residue system modulo $q$. Also note that
$$
K(u,v_i;p_i^{m_i}) = K\(u_i, v_i;p_i^{m_i}\),
$$
where the $v_i$ are determined by $v$ as before.
Then one can easily verify the identity
$$
\begin{aligned}
S_n(q)
 & = \sum_{u=1}^q \prod_{i=1}^s K(u,v_i;p_i^{m_i})^n \\
 & = \prod_{i=1}^s\sum_{u_i=1}^{p_i^{m_i}} K(u_i,v_i;p_i^{m_i})^n \\
 & =\prod_{i=1}^s S_n(p_i^{m_i}).
\end{aligned}
$$

So it suffices to study the power moments for Kloosterman sums for prime power modulus $q=p^r$. The case when $q=p$ is a prime (and thus $r=1$) has been studied extensively
in the literature from different points of views, see Robba \cite{Robba}, Katz \cite{Katz}, Evans \cite{Evans}, and Fu-Wan \cite{FuWan05}\cite{FuWan10} and the references listed there.
The goal of this paper is to try to understand $S_n(p^r)$ when $r>1$. Our main result is the following explicit formula for $S_n(p^r)$ if $r$ is suitably large.

\begin{theorem} Let $n\geq 1$ be odd. If $p$ is odd and $r  >  \log_p n +1$, then
$S_n(p^r)=0$. If $p=2$, then $S_n(2^r)=0$ for all $r >1$.
\end{theorem}


\begin{theorem} Let $n\geq 2 $ be even.
If $p$ is odd and  $r>  \log_p \frac{n}{2} +1$, then
$$
S_n(p^r) = \binom {n-1} {\frac n2 -1} \frac{p-1}p p^{(\frac n2 +1) r}.
$$
If $p=2$ and $r> \log_2 n + 2$,  we have
$$
S_n(2^r)= \binom {n-1} {\frac n2 -1}  2^{\frac n2 -2} 2^{(\frac n2 +1)r}.
$$
\end{theorem}

The bounds on $r$ in the conditions of the theorems are optimal. Indeed, we observed in several experiments that the clean formulas above for $S_n(p^r)$ do not hold when the lower bound on $r$ is not satisfied. When $n$ gets large compared to $p^r$, the behaviour of the moments becomes clearly more complicated, as illustrated by the following expressions, that we derived for $S_n(p^2)$. The \lq correction terms\rq\ depend on how $n$ behaves$ \mod p$.
(When $n$ is even, we do not provide an explicit proof, but it is similar to the proof when $n$ is odd.)

\begin{proposition}
Let $n$ and $p$ be odd and $p^2>n$. Then
$$
S_n(p^2) = -p^{n+1}\sum_{0\leq i\leq n-2, \ 2i \in \{ n-2, n-4\}\!\mod p} {n-2 \choose i}.
$$
Let $n$ be even and $p >\max \{2, \sqrt{\frac{n}{2}}\} $. Then
$$
S_n(p^2)= \binom {n-1} {\frac n2 -1} \frac{p-1}p p^{n+2}
 -p^{n+1} \sum^*_{0\leq i\leq n-2, \ i \in \{ \frac{n}{2}-1, \frac{n}{2}-2\}\!\mod p} {n-2 \choose i},
$$
where $\sum^*$ means the two obvious terms $i= \frac{n}{2}-1$  and $\frac{n}{2}-2$ are excluded.
\end{proposition}

\smallskip
The material is organized as follows.
In Section $2$,  we relate the power moment of Kloosterman sums for prime-power modulus ($q=p^r$) to the number of solutions of a certain equation
over $(\mathbb{Z}/p^r\mathbb{Z})^*$. The latter can be naturally studied in the framework of Igusa's zeta functions. When $n$ is odd, the associated hypersurface is fortunately
non-singular in characteristic zero (although it may be singular when reduced modulo $p^r$). In this case a generalized Hensel lemma can be used to show that the sequence $S_n(p^r)$ ($r=1, 2, \cdots$)  stabilizes to zero when $r$ is larger than an explicit
constant depending on $p$ and $n$.  To obtain the optimal condition in the theorem, we have to work a little harder and resort to a more delicate analysis.
When $n$ is even, the hypersurface is unfortunately singular in characteristic zero and the problem becomes significantly deeper.
We have to use an elaborate calculation of the Igusa zeta function via an explicit embedded resolution of singularities of the hypersurface.
In Section 3, we establish the link with the Igusa zeta function. In Section $4$, we work out the detailed calculation in the case $p$ is odd.
In Section $5$, we deal with the case $p=2$.
To our pleasant surprise,
the final results turn out to be quite nice and we always get a simple explicit formula for $S_n(p^r)$ when $r$ is suitably large.

\medskip
\section{Power moments for prime-power moduli}

\medskip
Let $q=p^r$, $p$ prime, $r\geq 2$ a positive integer. Let $\psi_q:\Z/q\Z\to\C$ be the additive character $x\mapsto e^{2\pi ix/q}$. We denote the $n$-th power moment of the classical Kloosterman sums modulo $q$ by
$$
S_n(q)=\sum_{\lambda\in\Z/q\Z}\(\sum_{x\in (\Z/q\Z)^*}\psi_q(x+\lambda/x)\)^n,\quad n\in\Z^+.
$$
Expanding the inner power and using the orthogonal property of additive characters, we have

\begin{eqnarray*}
S_n(q)
& = & \sum_{\lambda\in\Z/q\Z}\sum_{x_1,\ldots,x_n\in (\Z/q\Z)^*}
\psi_q\(x_1+\cdots+x_n+\lambda\(\frac{1}{x_1}+\cdots+\frac{1}{x_n}\)\) \\
& = & q\sum_{\frac{1}{x_1}+\cdots+\frac{1}{x_n}\equiv\, 0~({\rm mod}~q)}\psi_q(x_1+\cdots+x_n) \\
& = & q\sum_{\frac{1}{x_1}+\cdots+\frac{1}{x_{n-1}}+1\equiv\, 0~({\rm mod}~q)}
\sum_{x_n\in (\Z/q\Z)^*}\psi_q(x_n(x_1+\cdots+x_{n-1}+1)) \\
& = & q\sum_{\frac{1}{x_1}+\cdots+\frac{1}{x_{n-1}}+1\equiv\, 0~({\rm mod}~q)}
\(\sum_{x_n\in\Z/q\Z}-\sum_{x_n\equiv\,0~({\rm mod}~p)}\)\psi_q(x_n(x_1+\cdots+x_{n-1}+1)) \\
& = & q^2 \sum_{\substack{\frac{1}{x_1}+\cdots+\frac{1}{x_{n-1}}+1\equiv\,0~({\rm mod}~q)\\ x_1+\cdots+x_{n-1}+1\equiv\,0~({\rm mod}~q)}} 1 \\
&   & -~q\sum_{\frac{1}{x_1}+\cdots+\frac{1}{x_{n-1}}+1\equiv\,0~({\rm mod}~q)}
\sum_{x_n\in\Z/(q/p)\Z}\psi_{q/p}(x_n(x_1+\cdots+x_{n-1}+1)) \\
& = & q^2W_n(q) - q\cdot\frac{q}{p}\sum_{\substack{\frac{1}{x_1}+\cdots+\frac{1}{x_{n-1}}+1
\equiv\,0~({\rm mod}~q) \\ x_1+\cdots+x_{n-1}+1\equiv\,0~({\rm mod}~q/p)}}1 \\
& = & q^2W_n(q) - q^2/p\cdot (*).
\end{eqnarray*}
Here and in the sequel
$$W_n(p^k)= \sum_{\substack{\frac{1}{x_1}+\cdots+\frac{1}{x_{n-1}}+1\equiv\,0~({\rm mod}~p^k)\\ x_1+\cdots+x_{n-1}+1\equiv\,0~({\rm mod}~p^k)}}1$$
for $k \geq 1$, and $(*)$ is the similar sum on the last but one line above.
We first relate $(*)$ to $W_n(q/p)$. If $x_1,\ldots,x_{n-1}$ satisfies the system of congruences
$$
\begin{cases}
 x_1+\cdots+x_{n-1}+1\equiv 0\pmod{q/p} \\
 \frac{1}{x_1}+\cdots+\frac{1}{x_{n-1}}+1\equiv 0\pmod{q},
\end{cases}
$$
and we write $x_i\to x_i-\frac{q}{p}y_i$, $y_i\in\F_p$, then the second mod $q$ congruence  could be written as
$$
\frac{1}{x_1-\frac{q}{p}y_1}+\cdots+\frac{1}{x_{n-1}-\frac{q}{p}y_{n-1}}+1\equiv 0\pmod{q}.
$$
Noting that $(q/p)^2\equiv 0\pmod q$ for $q\neq p$, we have
$$
\frac{1}{x_1}+\cdots+\frac{1}{x_{n-1}}+1+\frac{q}{p}\(\frac{y_1}{x_1^2}+\cdots+\frac{y_{n-1}}{x_{n-1}^2}\)\equiv 0\pmod{q},
$$
implying
$$
\frac{p}{q}\(\frac{1}{x_1}+\cdots+\frac{1}{x_{n-1}}+1\)+\(\frac{y_1}{x_1^2}+\cdots+\frac{y_{n-1}}{x_{n-1}^2}\)\equiv 0\pmod{p}.
$$
Thus for fixed $x_1,\ldots,x_{n-1}$, the number of solutions in $(y_1,\ldots,y_{n-1})\in\F_p^{n-1}$ equals $p^{n-2}$. That is,
$$
(*)=\sum_{\substack{\frac{1}{x_1}+\cdots+\frac{1}{x_{n-1}}+1\equiv\,0~({\rm mod}~q)\\ x_1+\cdots+x_{n-1}+1\equiv\,0~({\rm mod}~q/p)}}1
=W_n(q/p)\cdot p^{n-2}.
$$
Similarly one easily computes $S_n(p)$ in terms of $W_n(p)$.
We summarize.

\begin{theorem}
If $q=p^r$ with $r\geq 2$, then
$$
S_n(q)=q^2W_n(q)-q^2p^{n-3}W_n(q/p)=q^2\(W_n(q)-p^{n-3}W_n(q/p)\).
$$
If $q=p$, then 
$$
S_n(p)=p^2W_n(p)-((p-1)^{n-1}+(-1)^n).
$$
\end{theorem}

\smallskip
For small $n$,  explicit formulas for $S_n(q)$ can be derived directly.
For $n=1$, it is clear that $S_1(q)=0$ for all $q$. For $n=2$, one checks that $W_n(q)=1$, and thus
$$
S_2(q) = q^2(1 -\frac{1}{p}).
$$
For $n=3$, one checks that $W_3(q)= 1 + (\frac{p}{3})$ if $p>2$ (and zero if $p=2$). This gives for $r\ge 2$ and all $p$ that
$$
S_3(q) = 0.
$$
Thus, we shall assume that $n\ge 4$ below.  (The case $n=4$ may also be do-able directly. At any rate, it follows from the explicit formulas in later sections).

Solving $x_{n-1}$ from the congruence $x_1+ \cdots + x_{n-1} + 1 \equiv 0 \pmod{q}$ and substituting it into the
congruence
$$
\frac{1}{x_1}+\cdots +\frac{1}{x_{n-1}} + 1 \equiv 0 \pmod{q},
$$
one finds that $W_n(q)$ is the number of $(\Z/q\Z)^*$-solutions of the equation
$$
(1+ x_1+x_2 +\cdots + x_{n-2})(1 +\frac{1}{x_1}+\frac{1}{x_2}+\cdots +\frac{1}{x_{n-2}}) = 1.
$$

Let $V_n(q)$  denote the number of $(\Z/q\Z)^*$-solutions of the zero set of the Laurent polynomial
\begin{equation}\label{torichypersurface}
g(x_1,\ldots, x_{n-1})=(x_1+x_2 +\cdots + x_{n-1})(\frac{1}{x_1}+\frac{1}{x_2}+\cdots +\frac{1}{x_{n-1}}) -1.
\end{equation}
Replacing $x_i$ by $x_ix_{n-1}$ in $g(x_1,\ldots, x_{n-1})$ for $1\le i\le n-2$, one checks that
$$
W_n(q) = \frac{1}{\phi(q)} V_n(q),
$$
where $\phi(q)=\phi(p^r) = p^{r-1}(p-1)$ is the Euler function. With this new notation, the first part of the previous theorem can be restated as follows.

\begin{theorem}\label{Sdescription}
If $q=p^r$ with $r\ge 2$, then
$$
S_n(q) = \frac{q^2}{\phi(q)}(V_n(q)-p^{n-2} V_n(q/p)).
$$
\end{theorem}

\smallskip
We now consider lifting solutions mod $q/p$ to solutions mod $q$, and then relate $V_n(q/p)$ to $V_n(q)$.
If $g(x_1,\ldots, x_{n-1})$ has no singular toric solution modulo $p$, then the Hensel lemma gives the recursive formula $V_n(q) = p^{n-2}V_n(q/p)$ and thus $S_n(q)=0$. We now check when the Laurent polynomial $g$ has no singular toric solution modulo $p$.

More generally, for a positive integer $k$, let $x=(x_1,\ldots, x_{n-1})$ be a singular toric solution of $g=0$
modulo $p^k$. That is,
$$
g(x_1,\ldots, x_{n-1})\equiv\frac{\partial g}{\partial x_1}(x)\equiv\cdots\equiv
\frac{\partial g}{\partial x_{n-1}}(x) \equiv 0 \pmod {p^k}.
$$
For $1\leq i\le n-1$, we deduce that
$$
\frac{\partial g}{\partial x_i}(x) = (\frac{1}{x_1}+\cdots +\frac{1}{x_{n-1}})
- \frac{x_1+\cdots +x_{n-1}}{x_i^2} \equiv 0 \pmod{p^k}.
$$
It follows that
$$
x_1^2 \equiv x_2^2 \equiv \cdots \equiv x_{n-1}^2 \pmod {p^k}.
$$
For $k\leq 2$ or $p>2$, this implies that $x_i\equiv\pm x_0\pmod {p^k}$ for some $x_0\in (\mathbb{Z}/p^k\mathbb{Z})^*$.
For $p=2$ and $k\geq 3$, we have the slightly weaker congruence $x_i\equiv\pm x_0\pmod {2^{k-1}}$ for some $x_0\in (\mathbb{Z}/2^{k-1}\mathbb{Z})^*$.
Let
$$
m_+ = \#\{1\le i\le n-1\,|\,x_i \equiv x_0 \pmod{p^k}\},
$$
$$
m_- = \#\{1\le i\le n-1\,|\,x_i\equiv -x_0 \pmod{p^k}\}
$$
if $k\leq 2$ or $p>2$, and analogously mod $2^{k-1}$ if $p=2$ and $k\geq 3$.
Then, at this singular point, for $p>2$ or $k\leq 2$, we have
$$
g(x_1,\ldots, x_{n-1}) \equiv (m_+ - m_-)^2 -1 \equiv 0 \pmod {p^k},
$$
and for $p=2$ and $k\geq 3$, we have
$$
g(x_1,\ldots, x_{n-1}) \equiv (m_+ - m_-)^2 -1 \equiv 0 \pmod {2^{k-1}}.
$$
Since $m_+ + m_- =n-1$, it follows that for $p=2$ and $k\geq 3$, we have
$$
(2m_+ -(n-1))^2 -1 = (2m_+ -n)(2m_+ -(n-2)) \equiv 0 \pmod {2^{k-1}}.
$$
For $p>2$ or $k\leq 2$, we have
$$
(2m_+ -(n-1))^2 -1 = (2m_+ -n)(2m_+ -(n-2)) \equiv 0 \pmod {p^k}.
$$
This shows that if $p^k > n$ (with $p>2$ or $k\leq 2$), $g$ has no singular toric solutions modulo $p^k$.
In the case $k=1$ and $n$ odd, this is impossible if either $p > n$ or $p=2$. We obtain the following.

\begin{theorem}\label{thm:SpecialCase}
Let $n$ be odd. Assume that $p>n$ or $p=2$. If $r\geq 2$, then
$$
V_n(p^r) = p^{n-2}\cdot V_n(p^{r-1}),\qquad S_n(p^r)=0.
$$
\end{theorem}

\medskip
\begin{question} Is there a direct proof of the statement $S_n(p^r)=0$ in Theorem~\ref{thm:SpecialCase}? When will we have $K(\lambda_1) = -K(\lambda_2)$?
\end{question}

This theorem shows that for odd $n$, the case $p=2$ is completely settled. We now assume that $n$ is odd and $p$ is also odd in the rest of this section.
We show that the above arguments can be refined to settle the more general case  for odd $n$ and $p^2>n$.  For this, we first show that the singular points in $\mathbb{F}_p$ never lift to points modulo $p^2$ and hence never lift to points modulo $p^r$ with $r\geq 2$. Let $(x_1,\ldots, x_{n-1})$ be a singular point modulo $p$. As above, we can take $x_i \equiv \pm x_0 \pmod p$. This singular solution modulo $p$ lifts to a solution modulo $p^2$  only if
$$
g(x_1,\ldots, x_{n-1}) \equiv (2m_+ -n)(2m_+ -(n-2)) \equiv 0 \pmod {p^2}.
$$
This is not possible as $n$ is odd, $0\leq m_+\leq n-1$ and $p^2>n$.
This implies that for $r\geq 3$, we still have
$$
V_n(p^r) = p^{n-2} V_n(p^{r-1}),\qquad S_n(p^r) =0.
$$
For $r=2$, let $N(n,p)$ denote the number of singular solutions modulo $p$. Taking $x_0=x_{n-1}$, one finds that
$$
N(n,p) = (p-1)\sum_{\substack{0\leq i\leq n-2\\ i \in\{\frac{n}{2}-1, \frac{n}{2}-2\}~{\rm mod}~p}} {n-2\choose i},
$$
where $i$ corresponds to $m_+-1$. We deduce that
$$
V_n(p^2) = p^{n-2}(V_n(p) - N(n,p)),\qquad S_n(p^2) = -\frac{p^{n+1}}{p-1}N(n,p).
$$
We summarize.

\begin{theorem}\label{thm:GeneralCase}
Let $pn$ be odd and $p^2>n$. If $r\geq 3$, then
$S_n(p^r)=0$.
If $r=2$, then
$$
S_n(p^2) = -p^{n+1}\sum_{0\leq i\leq n-2, \ i \in \{ \frac{n}{2}-1, \frac{n}{2}-2\} \mod p} {n-2 \choose i}.
$$
\end{theorem}

Note that the last sum is zero in the case $p>n$ and $n$ odd, consistent with the previous theorem.
The first part of the theorem can be further improved as follows.

\begin{theorem}\label{thm:GeneralCase}
If $np$ is odd and $r> \log_p n +1$, then
$
S_n(p^r)=0$.
\end{theorem}
\smallskip

\noindent
{\sl Proof.} Our assumption implies that $p^{r-1}>n$. We have shown  that for odd $pn$ with $p^{r-1}>n$,
the Laurent polynomial $g=0$ has no singular solutions modulo $p^{r-1}$. That is, there are no integers $x_i$ prime to $p$ such that $x=(x_1,\ldots, x_{n-1})$ satisfies
$$
g(x)\equiv\frac{\partial g}{\partial x_1}(x)\equiv\cdots\equiv\frac{\partial g}{\partial x_{n-1}}(x) \equiv 0 \pmod {p^{r-1}}.
$$
This means that any solution $x=(x_1,\ldots, x_{n-1})$ counted in $V_n(p^{r-1})$ must satisfy the inequality
$$
{\rm ord}_p\{\frac{\partial g}{\partial x_1}(x), \ldots, \frac{\partial g}{\partial x_{n-1}}(x)\} \le r-2.
$$

We claim that any solution $x=(x_1,\ldots, x_{n-1})$ counted in $V_n(p^{r-1})$ satisfies the stronger inequality
$$
k_x:={\rm ord}_p\{\frac{\partial g}{\partial x_1}(x), \ldots, \frac{\partial g}{\partial x_{n-1}}(x)\} < \frac{r-1}{2}.
$$
Otherwise, let $x=(x_1,\ldots, x_{n-1})$ be a solution counted in $V_n(p^{r-1})$ satisfying
$$
0< \frac{r-1}{2}  \leq k_x \le r-2.
$$
Let $y_i = x_i +p^{k_x} z_i$, where $z_i \in \mathbb{Z}$. Since $2k_x \geq  r-1$ and
$k_x={\rm ord}_p\{\frac{\partial g}{\partial x_1}(x), \ldots, \frac{\partial g}{\partial x_{n-1}}(x)\}$, the Taylor expansion shows that
$$g(y_1, \cdots, y_{n-1}) \equiv g(x_1,\cdots, x_{n-1}) \equiv 0 \pmod {p^{r-1}}.$$
Now, reducing $x$ modulo $p^{k_x}$, we see that $\{y_1, \cdots, y_{n-1}\}$ is a singular solution modulo $p^{k_x}$. One has as before the congruence
$$
y_1^2 \equiv y_2^2 \equiv \cdots \equiv y_{n-1}^2 \pmod {p^{k_x}}.
$$
This implies that $y_i\equiv\pm y_0\pmod {p^{k_x}}$ for some $0\leq y_0< p^{k_x}$.  We choose $y_i$ such that
$y_i = \pm y_0$ for all $1\leq i\leq n-1$.
Let
$$
m'_+ = \#\{1\le i\le n-1\,|\,y_i = y_0 \}, \
m'_- = \#\{1\le i\le n-1\,|\,y_i = -y_0 \}.
$$
Then
$$
0 \equiv g(y_1,\ldots, y_{n-1}) = (m'_+ - m'_-)^2 -1\pmod {p^{r-1}}. $$
This implies that
$$ (2m'_+ -n)(2m'_+ -(n-2)) \equiv 0  \pmod {p^{r-1}}.
$$
It contradicts our assumption that $p^{r-1} > n$ and $n$ is odd. The claim is proved.

Let
$$k_0= \max_x k_x=\max_x {\rm ord}_p\{\frac{\partial g}{\partial x_1}(x), \ldots, \frac{\partial g}{\partial x_{n-1}}(x)\},$$
where $x=(x_1,\dots, x_{n-1})$ runs over all solutions counted in $V_n(p^{r-1})$. The above claim shows that
$k_0 < (r-1)/2$, that is, $r \geq 2(k_0+1)$.
A more general Hensel lemma (see \cite{Segers}) implies that for $s \geq 2(k_0+1)-1$, we have
$$
V_n(p^s) = V_n(p^{2(k_0+1)-1})p^{(n-2)(s-2(k_0+1)+1)}.
$$
This is done by applying the general Hensel lemma only to those solutions modulo $p^{k_0+1}$ which can be lifted to
solutions modulo $p^s$.
Thus, for $s\geq 2(k_0+1)$,  we still have
$$
V_n(p^s) = p^{n-2}V_n(p^{s-1}).
$$
Since $r \geq 2(k_0+1)$, we can take $s=r$ and conclude that $S_n(p^r)=0$. \qed

\bigskip
When $n$ is even, it turns out that the Laurent polynomial $g$ always has singular toric solutions modulo $p^k$. This makes the determination of $S_n(p^r)$ via $V_n(p^r)$ much more difficult. In the last two sections we solve the problem using algebraic and geometric techniques from the study of Igusa zeta functions. First we explain the link with our problem in the next section.

\medskip
\section{Relation with Igusa zeta functions}

\medskip
In this section we describe the hypersurface in \eqref{torichypersurface} rather as the zero set of the {\em polynomial}
$$
h = \( (x_1 +x_2+\cdots + x_{n-1})(\frac{1}{x_1}+\frac{1}{x_2}+\cdots +\frac{1}{x_{n-1}}) - 1\)x_1x_2\cdots x_{n-1}.
$$

Classically one studies the behavior of the $V_n(p^r)$ through the properties of its generating series. We put
$$
P(t)=\(\frac{p-1}{p}\)^{n-1} + \sum_{r\geq 1} V_n(p^r) \(p^{-(n-1)}t\)^r \in \Q[[t]] ,
$$
where the constant $\(\frac{p-1}{p}\)^{n-1}$ and the factor $p^{-(n-1)}$ are the standard conventions, in order to relate $P(t)$ in a natural way with the Igusa zeta function of $h$. Igusa \cite{Igusa74} proved that $P(t)$ is in fact a rational function in $t$ through the study of that zeta function. We will obtain information about the poles of $P(t)$ by studying the Igusa zeta function of $h$, and then use the precise description of the $V_n(p^r)$ in terms of the poles (and their orders) of $P(t)$, as calculated in \cite{Segers}.

We introduce the version of the Igusa zeta function that we will use. We denote by $\Q_p$ and $\Z_p$ the field of $p$-adic numbers and the ring of $p$-adic integers, respectively, by $|\cdot|$ the standard $p$-adic norm on $\Q_p$ and by $dx$ the standard Haar measure on $\Q_p^k$. For $a\in\Z_p$ we denote by $\bar{a}$ its image in $\Z_p/p\Z_p \cong \F_p$, and similarly for $a\in \Z_p^k$ and $W\subset \Z_p^k$ we use the notation $\bar{a}$ and $\bar{W}$ for their images in $\F_p^k$.

\begin{definition} \rm
Let $f\in \Z_p[x_1,\dots,x_k]$ and let $W$ be a residual subset of $\Z_p^k$, that is, a disjoint union of residue classes$\mod p$. Then the Igusa zeta function associated to $f$ and $W$ is
$$
Z_W(f;s)=\int_W|f(x)|^s dx ,
$$
where $s\in \C$ with $\Re(s)>0$.
\end{definition}

Igusa \cite{Igusa74} showed in fact that $Z_W(f;s)$ is a rational function in $p^{-s}$, using an embedded resolution of singularities of $f$. Because of this result one considers $Z_W(f;s)$ as a function in $t=p^{-s}$ and writes $Z_W(f;t)$ for it. In fact $Z_W(f;t)$ contains the same information as the so-called Poincar\'e series
$$
P_W(f;t)= p^{-k}\# \bar{W} + \sum_{r\geq 1} V(f,W;p^r) \(p^{-k}t\)^r \in \Q[[t]] ,
$$
where $V(f,W;p^r)$ is the number of $k$-tuples $x \in
\(\Z_p/p^r\Z_p\)^k$ satisfying $f(x)\equiv 0 \pmod {p^r}$ and such
that the image of $x$ in $\F_p^k$ belongs to $\bar{W}$. Note that
the constant term is just the measure of $W$. More precisely one has
the relation
\begin{equation}\label{PversusZ}
P_W(f;t)=\frac{p^{-k}\# \bar{W} -tZ_W(f;t)}{1-t}
\end{equation}
by a straightforward adaptation of the proof of the standard case
\cite[Theorem 8.2.2]{Igusa00}. Note that, since $\Z_p/p^r\Z_p \cong
\Z/p^r\Z$, we have that $P(t)= P_{(\Z_p^*)^{n-1}}(h;t)$.

We now recall two techniques to compute the Igusa zeta function.

\bigskip
\noindent
{\sc The $p$-adic stationary phase formula.} We assume that at least one of the coefficients of $f$ does not belong to $p\Z_p$. (This can always be achieved by dividing $f$ by a suitable power of $p$.) Then we denote by $\bar{f}$ the non-zero polynomial over $\F_p$ obtained by reducing all the coefficients of $f$ modulo $p$.

Denote by $\bar{S}$ the subset of all $\bar{a}$ in $\bar{W}$ such that $\bar{f}(\bar{a})=0$ and $(\partial \bar{f}/\partial x_i)(\bar{a})=0$ for all $i\in \{1,\dots,k\}$, and by $S$ its preimage in $\Z_p^k$. Then \cite{Igusa94}\cite[Theorem 10.2.1]{Igusa00}
\begin{equation}\label{SPF}
Z_W(f;t)= p^{-k}(\# \bar{W} -N) + p^{-k}(N-\#\bar{S})\frac{(p-1)p^{-1}t}{1-p^{-1}t} + \int_S|f(x)|^s dx ,
\end{equation}
where $N$ is the number of zeroes of $\bar{f}$ in $\bar{W}$.

\bigskip
\noindent
{\sc Resolution of singularities.}  Let $\sigma:X\to \Q_p^k$ be an embedded resolution of singularities of $f$, where $X$ is a non-singular algebraic variety over $\Q_p$, $\sigma$ is a projective birational morphism, the inverse image of $\{f=0\}$ has simple normal crossings and $\sigma$ is an isomorphism outside that inverse image. Thus the irreducible components $E_i,i\in I,$ of $\sigma^{-1}\{f=0\}$ are nonsingular hypersurfaces, intersecting transversely. Note that at most $k$ different components $E_i$ contain a given point of $X$.  For $i\in I$ we denote by $N_i$ and $\nu_i -1$ the multiplicities of $E_i$ in the divisor of $\sigma^*f$ and of $\sigma^*(dx_1\wedge\dots\wedge dx_k)$, respectively. Then $Z_W(f;t)$ can be written as a rational function in $t$ with denominator $\prod_{i\in I}(1-p^{-\nu_i}t^{N_i})$, see \cite[Theorem 8.2.1]{Igusa00}. More precisely, $Z_W(f;t)$ is a sum of rational functions with denominator $\prod_{i\in J}(1-p^{-\nu_i}t^{N_i})$, where the $E_i, i\in J,$ have a nonempty intersection (and hence $\# J \leq k$).
Note that, by \eqref{PversusZ}, the same is then true for $P_W(f;t)$.

\medskip
There is an explicit formula of Denef \cite[Theorem 3.1]{Denef}, when a certain technical condition concerning the resolution $\sigma$ is satisfied. For the following notions we refer to \cite{Denef} for more information. To an algebraic set $V$ over $\Q_p$ is associated its reduction$\mod p$, being an algebraic set over $\F_p$ and denoted by $\bar{V}$. Also, to the map $\sigma$ one associates its reduction$\mod p$, being a morphism $\bar{\sigma}:\bar{X} \to \F_p^k$. When the restriction of $\sigma$ to $\sigma^{-1}W$ has good reduction$\mod p$ (see \cite{Denef} for this notion), we have
\begin{equation}\label{formula}
Z_W(f;t)= p^{-k}\sum_{J\subset I} c_J \prod_{i\in J}  \frac{(p-1)p^{-\nu_i}t^{N_i}}{1-p^{-\nu_i}t^{N_i}} ,
\end{equation}
where $c_J= \{\bar{x} \in \bar{X} \mid \bar{x} \in \bar{E_i} \text{ if and only if } i\in J, \text{ and } \bar{\sigma}(\bar{x}) \in \bar{W}\}$.
Here, to simplify notation,  we denote for a variety $\bar V$ over $\F_p$ the set of its $\F_p$-rational points by the same symbol $\bar V$.

\smallskip
In the next two sections we use these techniques to study $P(t)$ through the Igusa zeta
function associated to $h$.

\medskip
\section{Formula for $S_n(p^r)$ when $p$ is odd}

\medskip
{\em We assume in this section that $n$ is even and $p$ is odd.}
In fact we determined already when there exist $\bar{a}$ in $(\F_p^*)^{n-1}$ such that $\bar{h}(\bar{a})=0$ and $(\partial \bar{h}/\partial x_i)(\bar{a})=0$ for all $i\in \{1,\dots,n-1\}$. We use the notation $m_+$ and $m_-$ as before.
Replacing $x_0$ by $-x_0$ if necessary,  we may assume that $0\leq m_+ \leq n/2 -1$.
There exist such $\bar{a}$ in $(\F_p^*)^{n-1}$ if and only if
\begin{equation}\label{m+condition}
2m_+ \equiv n \pmod p \qquad\text{or} \qquad 2m_+ \equiv n-2 \pmod p .
\end{equation}

\medskip
Since $n$ is even,
\eqref{m+condition} is equivalent to $m_+ \equiv \frac n2 \pmod p$
or $m_+ \equiv \frac n2 -1 \pmod p$. When $p\geq \frac n2 +1$ this
happens if and only if $m_+= \frac n2 -1$. When $3\leq p \leq \frac
n2$ this happens for $m_+= \frac n2 -1$ and for at least one other
$m_+$, namely at least for $m_+= \frac n2 -p$.

\bigskip
We study the Igusa zeta function $Z(t)=Z_{(\Z_p^*)^{n-1}}(h;t)$. The
hypersurface $h=0$ in $(\Q_p^*)^{n-1}$ has singularities (of
multiplicity $2$) at $x_1^2=x_2^2=\cdots=x_{n-1}^2$. With a similar
argument as above these are the points where each $x_i=\pm x_0$ for
some $x_0\in \Q_p^*$ and $\#\{1\le i\le n-1\,|\,x_i=x_0\}=\frac n2 -1$.
Hence the singular locus of $h=0$ consists of $\binom {n-1}{\frac n2
-1}$ disjoint copies of $\Q_p^*$. One obtains an embedded resolution
$\sigma$ by blowing up with centres these lines; each exceptional
component $E_i$ is the product of such a centre $Z_i$ with a
$(n-3)$-dimensional projective space and has data $(N_i,
\nu_i)=(2,n-2)$. The strict transform $E_0$ of $\{h=0\}$ has data
$(N_0, \nu_0)=(1,1)$. We describe now the intersection of a fixed
$E_i$ with the strict transform. One easily computes that the
quadratic form
\begin{equation}\label{quadraticform}
q= \sum_{i=1}^{\frac n2 -1} x_i^2 + \sum^{n-2}_{\substack{ i,j=1 \\ i<j }} x_ix_j
\end{equation}
is the lowest degree term of a local equation of a transversal
section of the hypersurface $h=0$ at a singular point. Consequently
the intersection of $E_i$ with $E_0$ is the product of the centre
with the {\sl projective} variety determined by $q=0$.

Using Igusa's result above, we see that $Z(t)$ can be written as a rational function with denominator $(1-p^{-1}t)(1-p^{-n+2}t^{2})$.

\bigskip
\noindent
{\sc First case: $p\geq \frac n2 +1$.}

\medskip
One can check that $\sigma$ has good reduction$\mod p$, and hence we can apply Denef's formula \eqref{formula}.
Let $N$ denote the number of zeroes of $\bar{h}$  in $(\F_p^*)^{n-1}$. We claim that $p^{n-2}-1$ is the number of points of $\bar{E_i}$ mapping to $(\F_p^*)^{n-1}$ by $\bar{\sigma}$. Indeed, this is the product of $p-1$, being the number of points of $\bar{Z_i}\cap(\F_p^*)^{n-1}$, and the number of points of projective $(n-3)$-space over $\F_p$. Finally we denote by $Q$ the number of points of $\bar{E_i} \cap \bar{E_0}$ mapping to $(\F_p^*)^{n-1}$; it is the product of $p-1$ and the number of points on the {\sl projective} variety determined by $\bar{q}=0$. Then Denef's formula yields
\begin{equation}\label{h-formula}
\begin{aligned}
p^{n-1}Z(t)=& \ (p-1)^{n-1}-N + \big( N-\binom {n-1} {\frac n2 -1} (p-1)\big)\frac{(p-1)p^{-1}t}{1-p^{-1}t} \\
+&  \binom {n-1} {\frac n2 -1} ( p^{n-2} -1 - Q) \frac{(p-1)p^{-n+2}t^2}{1-p^{-n+2}t^2}\\
+& \binom {n-1} {\frac n2 -1}  Q \frac{(p-1)^2p^{-n+1}t^3}{(1-p^{-1}t)(1-p^{-n+2}t^2)} .
\end{aligned}
\end{equation}
More concretely, since $Q$ is also $1$ less than the number of points of the {\sl affine} variety determined by $\bar{q}=0$, we have by \cite[Theorem 9.2.1]{Igusa00} that
\begin{equation}\label{Q}
Q=(p^{\frac n2 -2}+1)(p^{\frac n2 -1}-1).
\end{equation}

\bigskip
\noindent
{\em General case: $n\geq 6$.}
It will turn out that we can write $Z(t)$, applying decomposition in partial fractions, in the form
\begin{equation}\label{ABC}
A + \frac B{1-p^{-1}t} + \frac C{1-p^{-\frac n2 +1}t}
\end{equation}
with $A,B,C$ constants. (Note that one expects a priori a term of the form $\frac {D+Et}{1-p^{-n +2}t^2}$. However, this term simplifies.)
A similar statement is then true for $P(t)$, yielding a concrete description of the behavior of $V_n(p^r)$ for $r\geq 1$.

We provide some details of this computation. Decomposing the last two terms of \eqref{h-formula} in partial fractions yields, as contribution to $\frac 1{1-p^{-n +2}t^2}$, the terms
$$\binom {n-1} {\frac n2 -1} ( p^{n-2} -1 - Q)(p-1) \quad\text{and}\quad -\binom {n-1} {\frac n2 -1}\frac {Q(p-1)^2}{p^{n-4}-1} (p^{n-4}+p^{-1}t),$$
respectively. Adding, dividing by $p^{n-1}$,  plugging in the
expression in \eqref{Q} for $Q$ and simplifying yields
\begin{equation}\label{simplified}
- \frac{\binom {n-1} {\frac n2 -1} (p-1)^2 (p^{\frac n2 -1} -1)}{p^{\frac n2 +1}(p^{\frac n2 -2} -1)} \cdot \frac{1+p^{-\frac n2 +1}t}{1-p^{-n +2}t^2} ,
\end{equation}
and indeed the last factor is equal to $\frac 1{1-p^{-\frac n2 +1}t}$.

In order to find the constant $C$ in the expression \eqref{ABC} for $P(t)$, we only need the similar constant in the expression for $Z(t)$. Using \eqref{PversusZ} one easily derives that $P(t)$ can be written in the form \eqref{ABC} with
\begin{equation}\label{C}
C= -\frac {\binom {n-1} {\frac n2 -1} (p-1)^2}{p^2 (p^{\frac n2 -2} -1)} .
\end{equation}
Looking at the main result in \cite{Segers} and its proof, we have for all $r\geq 1$ that
\begin{equation}\label{Segersformula}
V_n(p^r)= Bp^{(n-2)r} +  Cp^{\frac n2 r} .
\end{equation}
We compute by Theorem \ref{Sdescription} that
$$
\begin{aligned}
S_n(p^r)&=\frac{p^{2r+1}}{p-1} \( \frac{V_n(p^r)}{p^r} - p^{n-3}\frac{V_n(p^{r-1})}{p^{r-1}}  \) \\
&=\frac{p^{2r+1}}{p-1} \(Bp^{(n-3)r} +  Cp^{(\frac n2 -1) r} - p^{n-3} (Bp^{(n-3)(r-1)} +  Cp^{(\frac n2 -1)( r-1)})  \)  \\
&= \frac{p^{2r+1}}{p-1} C \(p^{(\frac n2 -1) r} - p^{\frac n2 -2 +(\frac n2 -1) r}  \) \\
&= C \frac {p(1-p^{\frac n2 -2})}{p-1} p^{(\frac n2 +1) r}
\end{aligned}
$$
for all $r\geq 2$. Plugging in \eqref{C} we obtain finally for all $r\geq 2$ that
\begin{equation}\label{Sformula}
S_n(p^r)= \binom {n-1} {\frac n2 -1} \frac{p-1}p p^{(\frac n2 +1) r}  .
\end{equation}

\bigskip
\noindent {\em Case $n=4$.}  In this special case a straightforward
calculation simplifies \eqref{h-formula} to
$$
Z(t)= \frac{p-1}{p^5}\cdot \frac {p^2(p^2-5p+7)+p(p^2-2p-5)t+(p^2+p+1)t^2}{(1-p^{-1}t)^2} ,
$$
yielding
$$
P(t)= \frac{p-1}{p^5}\cdot \frac {p^2(p^2-2p+1)+p(p^2-2p-2)t+(p^2+p+1)t^2}{(1-p^{-1}t)^2} .
$$
Decomposing $P(t)$ in partial fractions now results in the form
\begin{equation}\label{ABC2}
A + \frac B{1-p^{-1}t} + \frac C{(1-p^{-1}t)^2}
\end{equation}
with $A,B,C$ constants, and more precisely $C=3\frac{(p-1)^2}{p^2}$. In this case we have for all $r\geq 1$ by \cite{Segers} that
\begin{equation}\label{Segersformula2}
V_n(p^r)= \( (r+1)C + B\) p^{2r} .
\end{equation}
(Note that there is a typo in \cite{Segers} precisely at this point.
On the last line of page $4$ the numbers involving $e$ must be
augmented by $1$.) By Theorem \ref{Sdescription} we compute
$$
\begin{aligned}
S_n(p^r)&=\frac{p^{2r+1}}{p-1} \( \frac{V_n(p^r)}{p^r} - p\frac{V_n(p^{r-1})}{p^{r-1}}  \) \\
&=\frac{p^{2r+1}}{p-1} \( ((r+1)C+B) p^{r} - p (rC+B) p^{r-1}  \)  \\
&= \frac{p^{2r+1}}{p-1} C p^r \\
&= 3 \frac{p-1}p p^{3r}  .
\end{aligned}
$$
Note that this turns out to be exactly \eqref{Sformula} when substituting $n=4$.

\bigskip
\noindent
{\sc Second case: $3\leq p\leq \frac n2$.}

\medskip
(Hence $n\geq 6$.) We partition the $\bar{a}$ in $(\F_p^*)^{n-1}$ such that $\bar{h}(\bar{a})=0$ and $(\partial \bar{h}/\partial x_i)(\bar{a})=0$ for all $i\in \{1,\dots,n-1\}$ into the subsets $\bar{S_1}$, corresponding to $m_+ = \frac n2 -1$, and $\bar{S_2}$, corresponding to all other values of $m_+$. Let $S_1$ and $S_2$ denote their preimages in $\Z_p^{n-1}$, respectively.
The $p$-adic stationary phase formula \eqref{SPF} yields
\begin{equation}
\begin{aligned}
p^{n-1}Z(t)=& (p-1)^{n-1} -N + (N-\#(\bar{S_1}\cup\bar{S_2}))\frac{(p-1)p^{-1}t}{1-p^{-1}t} \\  &+p^{n-1} \int_{S_1}|h(x)|^s dx + p^{n-1} \int_{S_2}|h(x)|^s dx .
\end{aligned}
\end{equation}
In fact, the restriction of $\sigma$ to $S_1$ still has good
reduction$\mod p$, 
and by Denef's formula $p^{n-1}
\int_{S_1}|h(x)|^s dx$ equals the sum of the last two terms in
\eqref{h-formula}.

On the other hand, since $h=0$ has no singular points in $S_2$, we can write $p^{n-1} \int_{S_2}|h(x)|^s dx$ as a rational function in $t$ with denominator $1-p^{-1}t$. In general we cannot apply Denef's formula here; in particular we have no control over the degree of the numerator.
At any rate, decomposing $Z(t)$ and $P(t)$ in partial fractions, this time we can write $P(t)$ in the form
\begin{equation}
A'_n(t) + \frac {B'_n}{1-p^{-1}t} + \frac {C_n}{1-p^{-\frac n2 +1}t}  \qquad (n\geq 6),
\end{equation}
where $A'_n(t) \in \Q[t]$, $B'_n$ is a constant and $$C_n=-\frac {\binom {n-1} {\frac n2 -1} (p-1)^2}{p^2 (p^{\frac n2 -2} -1)}$$ as before.
By \cite{Segers} we still have similar expressions for $V_n(p^r)$ as in \eqref{Segersformula} and \eqref{Segersformula2}, but now they are only valid when $r$ is big enough, more precisely when $r>\deg A'_n(t)$.
We conclude that $V_n(p^r)$ is still given by the formula in \eqref{Sformula} when $r$ is big enough (with respect to $n$ and $p$).

We note that \eqref{Sformula} is also valid for $n=2$ and we summarize.

\begin{theorem} Let $n$ be an even positive integer.
Let $p$ be an odd prime number and $r\geq 2$.
If $p\geq \frac n2 +1$, then we have for all $r\geq 2$ that
$$
S_n(p^r)= \binom {n-1} {\frac n2 -1} \frac{p-1}p p^{(\frac n2 +1) r}.
$$
If $3\leq p\leq \frac n2$, the same formula is valid for $r$ big enough (depending on $n$ and $p$).
\end{theorem}

For $3\leq p\leq \frac n2$, the above theorem gives the precise information when $r$ is big enough. For small $r$,
the problem is caused by the integration over $S_2$, corresponding to points which are non-singular
over $\mathbb{Q}_p$, but become singular modulo $p$. This part can be handled
as in the second section when counting $V_n(p^r)$.  Combining the elementary method of that section
and the above Igusa zeta function calculation, we obtain the following additional results.

\begin{theorem} Let $n$ be an even positive integer and $p >\max \{2, \sqrt{\frac{n}{2}}\} $.
For $r \geq 3$, we have
$$
S_n(p^r)= \binom {n-1} {\frac n2 -1} \frac{p-1}p p^{(\frac n2 +1) r}.
$$
For $r=2$,  we have
$$
S_n(p^2)= \binom {n-1} {\frac n2 -1} \frac{p-1}p p^{n+2}
 -p^{n+1} \sum^*_{0\leq i\leq n-2, \ i \in \{ \frac{n}{2}-1, \frac{n}{2}-2\}  \mod p} {n-2 \choose i},
$$
where $\sum^*$ means the two obvious terms $i= \frac{n}{2}-1$  and $\frac{n}{2}-2$ are excluded.
\end{theorem}

Note that the second term is zero if $p \geq \frac{n}{2}+1$,
consistent with the previous theorem.

\begin{theorem} Let $n$ be an even positive integer and $p$ be odd.
For $r > \log_p \frac{n}{2} +1$, we have
$$
S_n(p^r) = \binom {n-1} {\frac n2 -1} \frac{p-1}p p^{(\frac n2 +1) r}.
$$
\end{theorem}

\medskip
Note that in the case $n$ even and $p$ odd, the modulo $p^k$ singularity condition
$$
(2m_+ -(n-1))^2 -1 = (2m_+ -n)(2m_+ -(n-2)) \equiv 0  \pmod {p^{k}}
$$
is equivalent to
$$
(m_+ -\frac{n}{2})(m_+ -(\frac{n}{2}-1)) \equiv 0  \pmod {p^{k}}.
$$
Thus, we can replace the previous condition $r-1 > \log_p n$ by the slightly weaker
condition  $r-1> \log_p \frac{n}{2} $.

\medskip
\section{Formula for $S_n(2^r)$}

\medskip
{\em We still assume in this section that $n$ is even.}
As before, the singular locus of the hypersurface
$h=0$ in $(\Q_2^*)^{n-1}$ consists of $\binom {n-1}{\frac n2
-1}$ disjoint copies of $\Q_2^*$, being
the points where each $x_i=\pm x_0$ for
some $x_0\in \Q_2^*$ and $\#\{1\le i\le n-1\,|\,x_i=x_0\}=\frac n2 -1$.
Blowing up with centres these lines yields an embedded resolution, and hence $Z(t)$ can be written as a rational function with denominator
$(1-2^{-1}t)(1-2^{-n+2}t^{2})$. But this resolution has bad reduction$\mod 2$.

\bigskip
\noindent
{\em General case: $n\geq 6$.}
We partition the integration domain $(1+2\Z_2)^{n-1}$ into (open and closed) pieces, where each piece contains at most one component of the singular locus. We can describe each such component $Z_J$ with its equations
$$
\begin{aligned}
&x_i=x_{n-1}  \qquad (i\in J) , \\
&x_i=-x_{n-1} \qquad (i \in \{1,\dots,n-2\}\setminus J)  ,
\end{aligned}
$$
where $J$ is a (uniquely determined) subset of $\{1,\dots,n-2\}$ with cardinality $\frac n2 -1$ or $\frac n2 - 2$.
We consider the neighborhood $U_J$ of $Z_J$ given by $x_i=x_{n-1} +4\Z_2$ $(i\in J)$ and $x_i=-x_{n-1} +4\Z_2$ $(i\notin J)$. Clearly $U_J$ and $U_{J'}$ are disjoint if $J\neq J'$. In order to compute
$$
\Int_J = \int_{U_J}|h(x)|^s dx
$$
we perform the (measure preserving) coordinate change $x_i=y_i+y_{n-1}$ $(i\in J)$, $x_i=y_i-y_{n-1}$ $(i\notin J)$, $x_{n-1}= y_{n-1}$, and we use the original description of the hypersurface. Then
$$
\begin{aligned}
\Int_J=\int_{(4\Z_2)^{n-2}\times (1+2\Z_2)} \big|&( \sum_{i=1}^{n-2} y_i \pm y_{n-1}) (\sum_{i\in J}\frac 1{y_i+y_{n-1}} \\&+ \sum_{i\notin J}\frac 1{y_i-y_{n-1}} +\frac 1{y_{n-1}})-1 \big|^s dy ,
\end{aligned}
$$
where in the first factor we have $+y_{n-1}$ (resp. $-y_{n-1}$) if $\#J=\frac n2 -1$ (resp. $\#J=\frac n2 -2$).
We further simplify the integral by \lq eliminating\rq\ the variable $y_{n-1}$. More precisely we perform the (also measure preserving) coordinate change $y_i=z_iz_{n-1}$ $(i\in \{1,\dots,n-2\}$, $y_{n-1}=z_{n-1}$, yielding
$$
\Int_J=\frac 12 \int_{(4\Z_2)^{n-2}} \big|( \sum_{i=1}^{n-2} z_i \pm 1) (\sum_{i\in J}\frac 1{z_i+1} + \sum_{i\notin J}\frac 1{z_i-1} + 1)-1 \big|^s dz ,
$$
where we used that $\int_{1+2\Z_2} dz_{n-1} = \frac 12$.
We can multiply the function within $|\cdot|$ with $\prod_{i\in J} (z_i + 1)\prod_{i\notin J} (z_i - 1)$ (having norm $1$ on the integration domain), in order to obtain a polynomial.
A straightforward computation yields that
$$
\( (\sum_{i=1}^{n-2} z_i \pm 1) (\sum_{i\in J}\frac 1{z_i+1} + \sum_{i\notin J}\frac 1{z_i-1} + 1)-1\) \prod_{i\in J} (z_i + 1)\prod_{i\notin J} (z_i - 1)
$$
is (up to sign) equal to
$$
2q(z) + g_{\geq 3}(z) ,
$$
where $g_{\geq 3}(z)$ contains only terms of degree at least $3$ and
\begin{equation}
\begin{aligned}
&q(z)= \sum_{i\notin J} z_i^2 + \sum_{1\leq i<j\leq n-2}z_i z_j \qquad \text{if } \# J=\frac n2 -1, \text{ and} \\
&q(z)= \sum_{i\in J} z_i^2 + \sum_{1\leq i<j\leq n-2}z_i z_j \qquad \text{if } \# J=\frac n2 -2  .
\end{aligned}
\end{equation}
Substituting $z_i=4x_i$ for $i=1,\dots,n-2$ yields
\begin{equation}\label{pieceofintegral}
\Int_J=\frac 12 \cdot \frac 1{4^{n-2}} \cdot 2^{-5s}\int_{(\Z_2)^{n-2}}\big| q(x) + 2 g'_{\geq 3}(x)\big|^s dx ,
\end{equation}
where $g'_{\geq 3}(x)$ contains only terms of degree at least $3$.

Note that the notation $q$ is consistent with (\ref{quadraticform}). In fact this last integrand has an isolated singularity in the origin, and blowing up at the origin yields an embedded resolution with good reduction$\mod 2$ and we can use Denef's formula. We can now proceed {\em completely analogously} as in the case $p\geq \frac n2 +1$. Comparing with the last two lines in (\ref{h-formula}), we claim that the contribution to $\int_{(\Z_2)^{n-2}}|q(x) + 2 g'_{\geq 3}(x)|^s dx$ involving $\frac 1{1-2^{-n+2}t^{2}}$ is
\begin{equation}\label{2contribution}
\frac 1{2^{n-2}}\(
  ( 2^{n-2} -1 - Q) \frac{2^{-n+2}t^2}{1-2^{-n+2}t^2}
+   Q \frac{2^{-n+1}t^3}{(1-2^{-1}t)(1-2^{-n+2}t^2)} \) ,
\end{equation}
where
$$
Q=(2^{\frac n2 -2}+1)(2^{\frac n2 -1}-1).
$$
In order to see this, we note the following.

($i$) The only difference is the factor $2^{n-2}$ (versus $p^{n-1}$). Indeed, now only $n-2$ variables are involved.

($ii$) For (\ref{h-formula}) the centres $Z_i$ were one-dimensional with $p-1$ as number of points of their reduction$\mod p$, and our present situation can be considered as a \lq transversal section\rq\ of the previous one. So we should a priori divide all \lq numbers of points\rq\ by $p-1$
to derive (\ref{2contribution}). But since here $p-1=1$ this makes no difference.

($iii$) The formula for $Q$ in \cite[Theorem 9.2.1]{Igusa00} is still valid for $p=2$ and for the two possible equations for $q$.

\smallskip
Arguing further as in the case $p\geq \frac n2 +1$, we see (compare with (\ref{simplified})) that the contribution of (\ref{2contribution}) to  $\frac 1{1-2^{-n+2}t^{2}}$ simplifies to
$$
- \frac{  (2^{\frac n2 -1} -1)}{2^{\frac n2}(2^{\frac n2 -2} -1)} \cdot \frac{1}{1-2^{-\frac n2 +1}t} .
$$
Combining this last expression with (\ref{pieceofintegral}), we see that the total contribution to $Z(t)$ involving $\frac{1}{1-2^{-\frac n2 +1}t}$ of all the integration domains $U_J$ is
$$
\begin{aligned}
\binom {n-1} {\frac n2 -1}\cdot \frac 12 \cdot \frac 1{4^{n-2}} \cdot t^{5} &\cdot \(- \frac{  (2^{\frac n2 -1} -1)}{2^{\frac n2}(2^{\frac n2 -2} -1)} \cdot \frac{1}{1-2^{-\frac n2 +1}t} \) \\
=& -\binom {n-1} {\frac n2 -1}\frac{  (2^{\frac n2 -1} -1)}{2^{5\frac n2 -3}(2^{\frac n2 -2} -1)}\cdot \frac{t^5}{1-2^{-\frac n2 +1}t} .
\end{aligned}
$$
Note that integrating $|h|^s$ over $(1+2\Z_2)^{n-1} \setminus \cup_J U_J$ will not contribute to a term involving $\frac{1}{1-2^{-\frac n2 +1}t}$ since $h=0$ is nonsingular there.

As in the previous cases our final aim is to determine $C_n$ in the description of $P(t)$ as
\begin{equation}\label{2ABC}
A_n(t) + \frac {B_n}{1-2^{-1}t} + \frac {C_n}{1-2^{-\frac n2 +1}t} ,
\end{equation}
where $A_n(t) \in \Q[t]$, and $B_n$ and $C_n$ are constants.
Since  $\frac{t^5}{1-2^{-\frac n2 +1}t}$ is the sum of a polynomial and $\frac{2^{5\frac n2 -5}}{1-2^{-\frac n2 +1}t}$, we conclude that, when writing $Z(t)$ in the form (\ref{2ABC}), the constant $C_n$ is equal to
$$
-\binom {n-1} {\frac n2 -1}\frac{  (2^{\frac n2 -1} -1)}{2^2(2^{\frac n2 -2} -1)}.
$$
Then, using as before (\ref{PversusZ}), one easily derives that $P(t)$ can be written in the form (\ref{2ABC}) with
$$
C_n= -\binom {n-1} {\frac n2 -1}\frac{  2^{\frac n2 -3}}{(2^{\frac n2 -2} -1)}.
$$
We conclude as before, by using Theorem \ref{Sdescription}, that
$$
\begin{aligned}
S_n(2^r)
& = C_n \cdot 2(1-2^{\frac n2 -2})2^{(\frac n2 +1)r} \\
& = \binom {n-1} {\frac n2 -1}  2^{\frac n2 -2} 2^{(\frac n2 +1)r}
\end{aligned}
$$
when $r$ is big enough (depending on $n$).

\bigskip
\noindent {\em Case $n=4$.}
Then the polynomial $h$ is simply $(x_1 + x_2)(x_1+x_3)(x_2+x_3)$ and one can compute in an elementary way that
$$
Z(t)=\frac{t^3(1-t+t^2)}{2^3(2-t)^2},
$$
and hence
$$
P(t)= \frac{4+t^2+t^3+t^5}{2^3(2-t)^2} .
$$
And then
$$
V_4(2^r) = \frac 32 (r-3) 2^{2r} \qquad\text{for } r>3
$$
and
$$
S_4(2^r) = 3\cdot 2^{3r} \qquad\text{for } r>4.
$$

We note again that this formula for $S_4(2^r)$ is compatible with the formula for $n\geq 6$, which is also compatible with
the formula for $n=2$ by the remark in section $4$. We summarize.

\begin{theorem} Let $n\geq 2$ be an even positive integer. Then
$$
S_n(2^r)= \binom {n-1} {\frac n2 -1}  2^{\frac n2 -2} 2^{(\frac n2 +1)r}
$$
when $r$ is big enough (depending on $n$).
\end{theorem}

Again, this result can be made more precise by using the elementary
method to remove the integration of $|h|^s$ over $(1+2\Z_2)^{n-1}
\setminus \cup_J U_J$. We can use the ideas in the proof of Theorem
\ref{thm:GeneralCase}, but in order to obtain an optimal bound, we
need more subtle arguments.

\begin{theorem} Let $n\geq 2$ be an even positive integer. For $r> \log_2 n + 2$,  we have
$$
S_n(2^r)= \binom {n-1} {\frac n2 -1}  2^{\frac n2 -2} 2^{(\frac n2 +1)r}.
$$
\end{theorem}

\smallskip
\noindent {\sl Proof.} Recall that for small $r$, the problem is
caused by points which are non-singular over $\Q_2$, but become
singular modulo $2$. Also, we saw that the singular locus of
 $g=0$ in $(\Q_2^*)^{n-1}$ consists of the points
 $(x_1,\ldots,x_{n-1})$
where each $x_i=\pm x_0$ for some $x_0\in \Q_2^*$ and $\#\{1\le i\le
n-1\,|\,x_i=x_0\}=\frac n2 -1$.

Consider odd integers $x_i$ such that
 $x=(x_1,\ldots, x_{n-1})$ satisfies $g(x) \equiv 0 \pmod
 {2^{r-1}}$, but such that $x \pmod 2$ does not lift to a singular
 solution over $\Z_2$.
 We claim that $x$ satisfies the inequality
$$
k_x:={\rm ord}_2\{\frac{\partial g}{\partial x_1}(x), \ldots,
\frac{\partial g}{\partial x_{n-1}}(x)\} < \frac{r}{2}.
$$
Otherwise, suppose that $$ 0< \frac{r}{2} \leq k_x .$$
 Let
$y_i = x_i +2^{\min(k_x,r-1)-1} z_i$, where $z_i \in \mathbb{Z}$. As before,
we want to argue using the Taylor expansion. In this case, an easily
verified but important fact is that all second partial derivatives
$\frac{\partial^2 g}{\partial x_i \partial x_j}(x)$ are congruent to
$0$ modulo $2$ (using only that the $x_i$ are odd).

Since $2\min(k_x,r-1) -1\geq r-1$ and $k_x={\rm ord}_2\{\frac{\partial
g}{\partial x_1}(x), \ldots, \frac{\partial g}{\partial
x_{n-1}}(x)\}$, the Taylor expansion shows that
$$g(y_1, \cdots, y_{n-1}) \equiv g(x_1,\cdots, x_{n-1}) \equiv 0 \pmod {2^{r-1}}.$$
Because $x$ is clearly a singular solution modulo $2^{\min(k_x, r-1)}$, one has
as in the fourth section the congruence
$$
x_1^2 \equiv x_2^2 \equiv \cdots \equiv x_{n-1}^2 \pmod {2^{\min(k_x, r-1)}},
$$
implying that $x_i\equiv\pm y_0\pmod {2^{\min(k_x, r-1)-1}}$ for some $y_0$ satisfying \linebreak

\noindent
$0\leq
y_0< 2^{\min(k_x, r-1)-1}$. We choose $y_i$ such that $y_i = \pm y_0$ for all
$1\leq i\leq n-1$. Let
$$
m'_+ = \#\{1\le i\le n-1\,|\,y_i = y_0 \}, \ m'_- = \#\{1\le i\le
n-1\,|\,y_i = -y_0 \},
$$
where $0\leq m'_+ \leq \frac n2 -1$. Then
$$
0 \equiv g(y_1,\ldots, y_{n-1}) = (m'_+ - m'_-)^2 -1\pmod {2^{r-1}}.
$$
 This implies that
$$ (2m'_+ -n)(2m'_+ -(n-2)) \equiv 0  \pmod {2^{r-1}},
$$
which is equivalent to
$$
(m'_+ -\frac{n}{2})(m'_+ -(\frac{n}{2}-1)) \equiv 0 \pmod {2^{r-3}}.
$$
Our assumption that $2^{r-2} > n$, or equivalently, $2^{r-3} > \frac
n2$, then implies that $m'_+ =\frac{n}{2}-1$. This contradicts the
condition that we imposed on $x \pmod 2$.
 The claim is
proved.

As in the proof of Theorem \ref{thm:GeneralCase}, we want to conclude using some Hensel lemma.
Let
$$k_0= \max_x k_x=\max_x {\rm ord}_2\{\frac{\partial g}{\partial x_1}(x), \ldots, \frac{\partial g}{\partial x_{n-1}}(x)\},$$
where $x=(x_1,\ldots, x_{n-1})$ runs over all solutions modulo $2^{r-1}$  such that $x \pmod 2$ does not lift to a singular
 solution over $\Z_2$. The above claim shows that $k_0 < r/2$.

If $k_0 < (r-1)/2$, the general Hensel lemma implies, as in the
proof of Theorem \ref{thm:GeneralCase}, that each such solution
modulo $2^{r-1}$ lifts to exactly $2^{n-2}$ solutions modulo $2^r$,
and so on. And then the contribution to $S_n(2^r)$ is zero.

When $k_0 = (r-1)/2$ (implying that $r$ is odd), we cannot invoke
the statement of the Hensel lemma, but in this case we can adapt its
classical proof with Taylor expansions, using again the crucial fact
that all second partial derivatives $\frac{\partial^2 g}{\partial
x_i \partial x_j}(x)$ are congruent to $0$ modulo $2$. More
precisely, let $x$ be a solution modulo $2^{r-1}$ as above with
moreover $k_x= (r-1)/2$. Let $y_i = x_i +2^{(r-1)/2} z_i$, where
$z_i \in \mathbb{Z}$. Looking at the Taylor expansion, requiring
that $g(y_1, \cdots, y_{n-1}) \equiv 0 \pmod {2^{r}}$ yields a
non-trivial linear relation modulo $2$ between $z_1,\ldots,z_{n-1}$,
that is, a non-trivial linear relation between their first digits.
Continuing this way we can still conclude that each such solution
$x$ modulo $2^{r-1}$ lifts to exactly $2^{n-2}$ solutions modulo
$2^r$, resulting again in a zero contribution to $S_n(2^r)$. (As
usual for the formal argument one has to start with solutions modulo
$2^{(r-1)/2}$ which can be lifted to solutions modulo $2^{r-1}$.)
  \qed


\end{document}